\documentclass[letterpaper,final,twoside,reqno,12pt]{amsart}
\usepackage{amsmath,amssymb}
\usepackage[all]{xy}
\usepackage{epic}
\usepackage{eepic}

\newcommand{\Z}{\ensuremath{\mathbb{Z}}}
\newcommand{\ZZ}{\ensuremath{\mathbb{Z} \oplus \mathbb{Z}}}
\newcommand{\Q}{\ensuremath{\mathbb{Q}}}
\newcommand{\R}{\ensuremath{\mathbb{R}}}
\newcommand{\RP}{\ensuremath{\mathbb{RP}}}
\newcommand{\MCG}{\ensuremath{\mathcal{MCG}}} 
\newcommand{\CC}{\ensuremath{\mathcal{C}}}

\newcommand{\GT}{\ensuremath{\mathcal{X}}}
\newcommand{\GD}{\ensuremath{\mathcal{D}}}
\newcommand{\GDG}{\ensuremath{\mathcal{D}_G}}

\newcommand{\LF}{\ensuremath{GVP(n)}} 
\newcommand{\LFnl}{\ensuremath{GVP_0(n)}} 

\newcommand{\tH}{\ensuremath{\widetilde{H}}} 
\newcommand{\tK}{\ensuremath{\widetilde{K}}} 

\newcommand{\co}{\colon\thinspace}

\newtheorem{theorem}{Theorem}[section]
\newtheorem{lemma}[theorem]{Lemma}

\newtheorem{corollary}[theorem]{Corollary}

\newtheorem{remark}[theorem]{Remark}

\DeclareMathOperator{\Aut}{Aut}
\DeclareMathOperator{\Isom}{Isom}
\DeclareMathOperator{\Out}{Out}

\DeclareMathOperator{\Comm}{Comm}

\DeclareMathOperator{\cd}{cd}
\DeclareMathOperator{\hi}{h}

\begin{document}


\title{A fixed point theorem for deformation spaces of $G$--trees}
\author{Matt Clay}
\address{Department of Mathematics\\ 
University of Utah\\ 
Salt Lake City, UT 84112-0090 USA}
\email{clay@math.utah.edu}
\date{\today}

\begin{abstract}
For a finitely generated free group $F_n$, of rank at least 2, any
finite subgroup of $\Out(F_n)$ can be realized as a group of
automorphisms of a graph with fundamental group $F_n$.  This result,
known as $\Out(F_n)$ realization, was proved by Zimmermann, Culler and
Khramtsov.  This theorem is comparable to Nielsen realization as
proved by Kerckhoff: for a closed surface with negative Euler
characteristic, any finite subgroup of the mapping class group can be
realized as a group of isometries of a hyperbolic surface.  Both of
these theorems have restatements in terms of fixed points of actions
on spaces naturally associated to $\Out(F_n)$ and the mapping class
group respectively.  For a nonnegative integer $n$ we define a class
of groups (\LF) and prove a similar statement for their outer
automorphism groups.
\end{abstract}

\maketitle


{AMS subject classification, primary: 20E08, secondary: 20F65, 20F28}

{keywords: $G$--tree, deformation space, $\Out(F_n)$ realization,
Nielsen realization}

\vskip 0.5cm Let $\Sigma$ be a closed surface with negative Euler
characteristic.  The mapping class group $\MCG (\Sigma)$ acts on
Teichm\"uller space $T_\Sigma$, the space of hyperbolic metrics on
$\Sigma.$ A stabilizer in this action is an isometry group of some
hyperbolic metric on $\Sigma$.  Such groups must be finite.  It is
theorem of Kerckoff \cite{ar:Ke}, known as Nielsen realization, that
any finite subgroup of $\MCG(\Sigma)$ can be realized as a group of
isometries for some hyperbolic metric on $\Sigma$.  Therefore the
finite subgroups of $\MCG (\Sigma)$ are exactly the subgroups with
fixed points in $T_\Sigma$.  In a similar manner, for a finitely
generated free group $F_n$ of rank $n \geq 2$, the outer automorphism
group $\Out(F_n)$ acts on Culler and Vogtmann's Outer space $X_n$
\cite{ar:CV}.  A stabilizer in this action are is an isometry group of
a some metric graph with fundamental group $F_n$.  It is a theorem of
Zimmermann \cite{ar:Z}, Culler \cite{ar:Cu} and Khramtsov
\cite{ar:Kh}, known as $\Out(F_n)$ realization, that any finite
subgroup of $\Out(F_n)$ can be realized as a group of isometries of
some metric graph with fundamental group $F_n$.  Thus as for
$\MCG(\Sigma)$ and $T_\Sigma$, the finite subgroups of $\Out(F_n)$ are
exactly those subgroups with a fixed point in $X_n$.

For a nonnegative integer $n$ we introduce a class of groups denoted
\LF, and prove a similar realization statement for their outer
automorphism groups.  In other words, for every group $G \in \LF,$
there is a naturally associated space on which $\Out(G)$ acts and we
are able to determine that certain subgroups of $\Out(G)$ related to
stabilizers have a fixed point.  For $n=0,1$ we show that any group
which is commensurable to a subgroup of a stabilizer actually fixes a
point (Corollary \ref{col:n=01}).  The class $GVP(0)$ is the class of
virtually finitely generated free groups of rank at least 2, thus our
result is a generalization of $\Out(F_n)$ realization.  In general, we
are only able to show that subgroups of $\Out(G)$ commensurable to
polycyclic subgroups of stabilizers fix a point.

We define $\LFnl$ as the class of groups which act on a locally finite
simplicial trees without an invariant point or line, such that the edge
stabilizers are virtually polycyclic subgroups of Hirsch length $n$.
The subset of groups in \LFnl \ where this action is irreducible and
cocompact is denoted \LF.  In the first section for any finitely
generated group $G$ we describe topological spaces \GD \ on which
certain subgroups of $\Out(G)$ act.  These spaces are contractible in
most cases.  In particular, for $G \in \LF$ we describe a contractible
topological space \GDG \ on which the full group $\Out(G)$ acts.  Our
main theorem regarding this action is analogous to Nielsen realization
and $\Out(F_n)$ realization.

\vskip 0.5cm
\noindent
\textbf{Main Theorem.} \textit{Suppose $G \in \LF$ and $K$ is a
polycyclic subgroup of $\Out(G)$ which fixes a point in \GDG.  If $H$
is a subgroup of $\Out(G)$ commensurable with $K$, then $H$ fixes a
point in \GDG.}

\vskip 0.5cm The proof of the above is similar to the proof for finite
subgroups of $\Out(F_n)$.  We review how to prove $\Out(F_n)$
realization.  Starting with a finite subgroup $K$ of $\Out(F_n),$ lift
this to the subgroup \tK \ in $\Aut(F_n)$.  Then \tK \ is virtually
free, hence \tK \ acts cocompactly on a simplicial tree $T$ with
finite stabilizers by Stallings' theorem \cite{ar:SW}.  This induces a
cocompact free action of $F_n \subseteq \tK$ on $T$.  Thus the finite
group $K = \tK / F_n$ acts on the quotient graph $T / F_n$, which
represents a point in Outer Space.  Hence this point is fixed by $K$.

We seek to mimic this proof.  The ingredient we will need is an analog
to Stallings' theorem, i.e. we need to know when can we raise a
splitting of a finite index subgroup to the whole group.  For the
special case we consider, this question has an answer due to Dunwoody
and Roller \cite{ar:DR}.  We then show that any group which contains a
finite index subgroup in \LF \ is in fact itself in \LF.  Finally, if
$G \in \LF$ and $K$ is a polycyclic subgroup of $\Out(G)$ which
stabilizes a point in \GDG, we show that \tK, the lift of $K$ to
$\Aut(G),$ is in $GVP(n')$ for some $n'$, inducing an action
of $G$.  Thus we can proceed as above for $\Out(F_n)$.

Originally, we were only concerned only with a proof of realization
for generalized Baumslag--Solitar (GBS) groups, the torsion-free groups in
$GVP(1)$.  However, in doing so it became necessary to prove
some statements in greater generality, which provided a proof for any
\LF --group.

\vskip 0.2 in
\noindent
\textbf{Acknowledgements}. I would like to thank my advisor Mladen
Bestvina for suggesting this line of research and for his
encouragement.  Thanks are also due to Michah Sageev for discussions
related to this work, to Gilbert Levitt for the correct references to
$\Out(F_n)$ realization and to the referee for shorting an early
version of this paper by pointing to Theorem \ref{th:fis}.


\section{Deformation Spaces of $G$--trees}\label{sc:ds}

For a finitely generated group $G$, a $G$--\textit{tree} is a metric
simplicial tree on which $G$ acts by isometries.  We say two
$G$--trees $T$ and $T'$ are equivalent if there is a $G$--equivariant
isometry between then.  When we speak of a $G$--tree we will always
mean the equivalence class of the $G$--tree.  A subgroup is an
\textit{elliptic} subgroup of $T$ if it fixes a point in $T$.  Given a
$G$--tree there are two moves one can perform to the tree that do not
change whether subgroups of $G$ are elliptic.  These moves correspond
to the isomorphism $A \cong A *_C C$ and are called \textit{collapse}
and \textit{expansion}.  For a detailed description of the moves see
\cite{ar:F1}.  In \cite{ar:F1} Forester proves the converse, namely if
two cocompact $G$--trees have the same elliptic subgroups, then there
is a finite sequence of collapses and expansions (called an
\textit{elementary deformation}) transforming one $G$--tree to the
other.

We let \GT \ denote a maximal set of $G$--trees which are related by an
elementary deformation.  By the theorem of Forester mentioned above,
an equivalent definition is as the set of all $G$--trees that have the
same elliptic subgroups as some fixed $G$--tree.  This set \GT \ is
called an \textit{unnormalized deformation space}.  We will always
assume that the $G$--trees are minimal and irreducible and that $G$
acts without inversions.

There is an action of $\R^+$ on \GT \ by scaling, the quotient is
called a \textit{deformation space} and denoted \GD.  We \cite{ar:C}
and independently Guirardel and Levitt \cite{ar:GL}, \cite{ar:GL2}
have shown that for a finitely generated group, if the actions in \GD
\ are irreducible and there is a reduced $G$--tree with finitely
generated vertex stabilizers, then \GD \ is contractible.  The
topology for the preceding statement is the axes topology induced
from the embedding $\GD \to \RP^\CC$ where $\CC$ is the set of all
conjugacy classes of elements in $G$, or equivalently the
Gromov--Hausdorff topology.  See \cite{ar:C} for details.

In general, the space \GD \ is acted on only by a subgroup of
$\Out(G)$, where the action is precomposition.  This subgroup is the
subgroup of $\Out(G)$ which permutes the conjugacy classes of elliptic
subgroups associated to \GD.

If $G \in \LF$ then there is a locally finite $G$--tree $T$ where all
of the stabilizers are virtually polycyclic subgroups of Hirsch length
$n$.  We will show in the next section (Lemma \ref{lm:lf-tree}) that
the set of elliptic subgroups for this action is invariant under all
automorphisms of $G$.  Hence the deformation space containing $T$ is
invariant under $\Out(G)$.  We denote this space as \GDG.  Notice that
if $G = F_n$, then $\GDG = X_n$, Culler and Vogtmann's Outer space.


\section{Virtually Polycyclic Groups and the Class of \LF --Groups}
\label{sc:lf-groups}

A group $G$ which admits a filtration $\{1 \} = G_0 \triangleleft G_1
\triangleleft \ldots \triangleleft G_n = G$ with $G_{i-1}
\triangleleft G_i$ normal and $G_i / G_{i-1}$ cyclic is called
\textit{polycyclic}.  The \textit{Hirsch length}, denoted $\hi G$, of
a polycyclic group $G$ is the number of infinite cyclic factors in the
above filtration.  This is an invariant of $G$.  If $G$ is polycyclic
and $H$ is any finite index subgroup then $\hi H = \hi G$.  In fact,
if $H$ is a subgroup of $G$, then $\hi H \leq \hi G$ with equality if
and only if $H$ has finite index in $G$.  This allows us to define the
Hirsch length of a group which contains a polycyclic group as a finite
index subgroup.  Such groups are called \textit{virtually polycyclic}.
These groups are also referred to as polycyclic-by-finite groups in
the literature.  Note that if $1 \to K \to G \to H \to 1$ is a short
exact sequence then $H$ and $K$ are virtually polycyclic if and only
if $G$ is.  In this case, the Hirsch lengths satisfy $\hi G = \hi H +
\hi K$.

As mentioned in the introduction, originally the main theorem was only
intended for generalized Baumslag--Solitar groups.  A group $G$ is a
\textit{generalized Baumslag--Solitar (GBS) group} if there is a
cocompact $G$--tree where the stabilizer of any point is isomorphic to
\Z.  As the only \Z \ subgroups of \Z \ are necessarily of finite
index, this $G$--tree must necessarily be locally finite.
Equivalently, $G$ is a GBS group if it admits a graph of groups
decomposition where all of the edge groups and vertex groups are
isomorphic to \Z.  Hence the nonelementary GBS groups (i.e. $G
\neq\Z,\ZZ$ or the Klein-bottle group) are the torsion-free groups in
$GVP(1)$.

Such groups were first studied by Kropholler \cite{ar:K}, where it is
shown that GBS groups are the only finitely generated groups of
cohomological dimension two that contain an infinite cyclic subgroup
which intersects each of its conjugates in a finite index subgroup.
It is clear that for a GBS group any vertex group in the above
mentioned graph of groups decomposition satisfies this condition.
Forester's Lemma 2.5 in \cite{ar:F2} (a generalization of which
appears as Lemma \ref{lm:lf-tree} below) implies that when the action
does not have an invariant line, the elliptic subgroups are the only
subgroups which satisfy this condition.  As this condition is
algebraic, the set of elliptic subgroups is invariant under all
automorphisms of $G,$ hence we can talk about an $\Out(G)$--invariant
deformation space.  We now generalize this fact to any \LF --group.

Recall that a group is called \textit{slender} if every subgroup
is finitely generated.  Virtually polycyclic groups are slender.
Slenderness of a group $G$ is equivalent to every subgroup $H
\subseteq G$ having property AR: whenever $H$ acts on a simplicial
tree, $H$ either stabilizes a point or has an axis \cite{ar:DSa}.
Throughout the following, we use the notation $H^g = gHg^{-1}$.  We
say two subgroups $H,H'$ of $G$ are \textit{commensurable} if $H \cap
H'$ has finite index in both $H$ and $H'$.  The commensurator of a
subgroup $H \subseteq G$ is $\Comm_G(H) = \{ g \in G \ | \ H \mbox{ is
commensurable with } H^g \}$.

\begin{lemma}\label{lm:lf-tree}{\rm (Forester \cite{ar:F2})}
Let $T$ be a locally finite $G$--tree such that the stabilizer of any
point in $T$ is slender.  If $T$ does not contain a $G$--invariant
line, then a subgroup $H \subseteq G$ is elliptic if and only if $H$
is contained in a subgroup $K$, where $K$ is slender and
$\Comm_G(K)=G$.
\end{lemma}

\begin{proof}
As $T$ is a locally finite simplicial tree, any vertex stabilizer is
commensurable to all of its conjugates.  Hence, if $H$ is elliptic, it
is contained in a vertex stabilizer $K$ which satisfies the conclusion
of the lemma.  

For the converse suppose $H \subseteq G$ does not act
elliptically and $H$ is contained in a slender subgroup $K$.  Hence
$K$ does not act elliptically.  $L_K$ be its axis.  Then the axis of
$K^g$ is $gL_K$.  If $K$ and $K^g$ are commensurable, then they have
the same axis so $L_K = gL_K$.  Hence if $K$ is commensurable to all
of its conjugates, then $L_K$ is a $G$--invariant line.
\end{proof}

Thus for such actions the elliptic subgroups are determined
algebraically.  In particular, the elliptic subgroups for these
actions are invariant under $\Aut(G)$.  When the action is cocompact,
we can talk about an $\Out(G)$--invariant deformation space, denote
this space \GDG. Hence every point in \GDG \ is a locally finite
$G$--tree where the stabilizers are virtually polycyclic of Hirsch
length $n$.  Since the $G$--trees in \GDG \ are locally finite, \GDG \
is a locally finite complex.  Hence all of the stabilizers are
commensurable.  Our realization statement (Main Theorem) is a partial
converse to this in the general case and the full converse if $n=0$ or
1.

For $G \in \LF$, as these actions are irreducible, $G$ contains a free
subgroup of rank 2.  Thus if $G$ acts on a tree $T$ with virtually
polycyclic stabilizers of Hirsch length $n$, then $T$ cannot be a
line.  This will be used without further mention.  We have another
lemma essentially due to Forester about the splittings of \LF --groups
as amalgams over virtually polycyclic groups $K$ with $\hi K =n$.  Say
that $G$ \textit{splits over $K$} if $G$ can either be written as a
nontrivial free product with amalgamation $G = A *_K B$ or as an
HNN-extension $G = A *_K$.

\begin{lemma}\label{lm:lf-ell-ell}{\rm (Forester \cite{ar:F2})}
Suppose $G \in \LF$ and $T \in \GDG$.  If $G$ splits over a virtually
polycyclic subgroup $K$ with $\hi K = n,$ then $K$ fixes a point in
$T$.  Moreover, the vertex group(s) in this splitting are finitely
presented and either \LFnl--groups or virtually polycyclic with Hirsch
length $n$ or $n+1$.
\end{lemma}

\begin{proof}
Let $Y$ be the Bass--Serre tree for the splitting of $G$ over $K$ and
$H$ a vertex stabilizer for the $G$--tree $T$.  Then similarly to
Lemma \ref{lm:lf-tree}, $H$ must act elliptically on $Y$ as $Y$ cannot
contain a $G$--invariant line.  Let $y \in Y$ be a vertex fixed by $H$
and $e$ an edge stabilized by $K$.  There is some $g \in G$ such that
$e$ separates $y$ from $gy$.  As $H$ and $H^g$ are commensurable,
there is a finite index subgroup $H' \subseteq H$ stabilizing $e,$
hence contained in $K$.  As both $K$ and $H'$ have Hirsch length $n,$
$H'$ has finite index in $K,$ hence $K$ an $H$ are commensurable.
Thus as $H$ fixes a point in any $G$--tree $T \in \GDG$, so does $K$.

As for the moreover, suppose $A$ is a vertex group for the splitting
of $G$ over $K$.  We examine how $A$ acts on $T$.  If $A$ fixes a
point, then $A$ is virtually polycyclic with $\hi A = n$.  If $A$ has
an invariant line on which it acts nontrivially, then there is a short
exact sequence $1 \to K' \to A \to \Z \to 1$ or $1 \to K' \to A \to
\Z_2 * \Z_2 \to 1$ where $K'$ is commensurable to $K$.  Hence $A$ is
virtually polycyclic and $\hi A = n+1$.  Otherwise this action implies
that $A \in \LFnl$.

To see that $A$ is finitely presented we can assume that $A \in
\LFnl$.  As $K$ acts elliptically in $T$, using the action of $A$ on
$T$ we can refine the splitting of $G$ over $K$ to get a graph of
groups decomposition for $G$ that includes the graph of groups
decomposition of $A$ with virtually polycyclic vertex and edge groups
of Hirsch length $n$.  As $G$ is finitely generated, after reducing we
can assume that the graph of groups decomposition for $G$, hence the
graph of groups decomposition of $A$, is a finite graph.  Thus $A$ can
be expressed as a finite graph of groups where all of the vertex and
edge groups are virtually polycyclic of Hirsch length $n$.  In
particular, $A$ is finitely presented.
\end{proof}

We record some properties about \LF --groups that will be used in section 
\ref{sc:rvts}.

\begin{lemma}\label{lm:prop}
Let $G \in \LF$ then:

\begin{itemize}

\item[1.] $\cd_{\Q} G = n+1$;

\item[2.] $G$ does not split over a virtually polycyclic group of
Hirsch length less than $n$; and

\item[3.] the center of $G$, $Z(G)$, is a virtually polycyclic subgroup
with $\hi Z(G) \leq n$.  The quotient $G / Z(G)$ is in
$GVP(n')$ for $n' = n - \hi Z(G)$.

\end{itemize}

\end{lemma}

\begin{proof}
For 1. and 2. see Bieri sections 6 and 7 \cite{ar:B}.

To see 3., let $T \in \GDG$.  As the action is irreducible, $Z(G)$
must act trivially on $T$ \cite{ar:BJ}.  Hence $Z(G)$ is a virtually
polycyclic subgroup with $\hi Z(G) \leq n$. Also, we have an induced
irreducible cocompact action of $G / Z(G)$ on $T$, where the
stabilizers are the quotients of the stabilizers for the $G$--action
by $Z(G)$.  Hence $G /Z(G)$ is in $GVP(n')$ where $n' = n - \hi Z(G)$.
\end{proof}

\begin{remark}\label{rm:prop}{\rm
As cohomological dimension is an invariant of the group, if $G \in
\LF,$ then $G \notin GVP(n')$ for $n \neq n'$. 
}
\end{remark}


\section{Promoting Finite Index Splittings}\label{sc:fis}

The main step in proving $\Out(F_n)$ realization is to use Stallings'
theorem to get a splitting of the virtually free subgroup of
$\Aut(F_n)$ which is the lift of some finite subgroup in $\Out(F_n)$.
In the present setting we will need an analog of Stallings' theorem to
tell us when a splitting of a finite index subgroup $H \subseteq G$
over $K \subseteq H$ implies a splitting of the whole group $G$ over
some subgroup $K'$ which is commensurable to $K$.  In general, we
cannot expect a splitting of $G$.  In our special case though, the
answer is given by the following theorem of Dunwoody and Roller
\cite{ar:DR} as stated by Scott and Swarup \cite{ar:SS}.  The ends of
the pair of groups $K \subseteq G$, denoted $e(G,K)$, is the number of
ends of $\Gamma/K$ where $\Gamma$ is a Cayley graph for $G$.  If $G$
splits over a subgroup $K$, then $e(G,K) > 1$.  See \cite{ar:SS} or
\cite{ar:SW} for these notions.

\begin{theorem}\label{th:fis}{\rm (Dunwoody--Roller \cite{ar:DR} as 
stated in \cite{ar:SS})} If $G,K$ are finitely generated subgroups
with $e(G,K)>1$ and if $\Comm_G(K) = G$, then $G$ splits over a
subgroup commensurable to $K$.
\end{theorem}

We can now prove our analog to Stallings' Theorem.

\begin{theorem}\label{th:st}
Let $G$ be a finitely presented group which has a finite index
subgroup $H \in \LF,$ then $G \in \LF$.
\end{theorem}

\begin{proof}
If $n=0$ then this is Stallings' theorem \cite{ar:SW}, so we assume
that $n \geq 1$.

Let $T \in \GD_H$ and $K$ be an edge stabilizer of $T$.  As finite
index subgroups of $H$ are in \LF, we can assume that $H$ is normal in
$G$.  By Lemma \ref{lm:lf-tree}, elliptic subgroups of $H$ are
invariant under automorphisms of $H$, hence $\Comm_G(K) = G$.  As
$e(H,K) > 1$ and $H$ is a finite index subgroup of $G$, we must have
that $e(G,K) > 1$.  Then by Theorem \ref{th:fis}, $G$ splits over a
subgroup $K'$ commensurable with $K$.  Let $T'$ be the Bass--Serre
tree for this splitting of $G$ over $K'$.  If $T'$ is locally finite
then we are done as the vertex and edge stabilizers for $G$ acting on
this tree are then commensurable to $K'$ hence virtually polycyclic
subgroups of Hirsch length $n$, therefore $G \in \LF$.

Suppose that $T'$ is not locally finite, we now show that we can split
a vertex group for the graph of group decomposition induced by $H$
acting on $T'$.  As $T'$ is not locally finite there is a vertex group
$H_v$ which is not a virtually polycyclic subgroup of Hirsch length at
most $n$.  Suppose $G_v$ is the vertex group under the $G$--action,
thus $H_v$ has finite index in $G_v$. As the induced action of $H_v$
on $T$ is nontrivial, we get a graph of groups decomposition for
$H_v$.  Then we can collapse this graph of groups decomposition to get
a splitting of $H_v$ over some edge stabilizer $K_v$.  Denote the
Bass--Serre tree for this splitting as $T_v$.  As $K_v$ is
commensurable to $K$ we have $\Comm_G(K_v) = G$, thus $\Comm_{G_v}(K_v) =
G_v$.  As $H_v$ is finitely generated, $G_v$ is also.  Therefore by as
above by Theorem \ref{th:fis}, $G_v$ splits as an amalgam over $K_v'$
which is commensurable to $K_v$ hence also $K'$.

As $K_v'$ and $K'$ are commensurable, in the Bass--Serre tree
associated to the splitting of $G_v$ over $K_v'$, $K'$
acts elliptically.  This allows us to refine the one edge splitting of
$G$ over $K'$ to get a two edge splitting of $G$ over $K'$ and $K_v'$.
Once again, we have a Bass--Serre tree $T_0$ associated to this graph
of groups decomposition, and the action of $H$ on $T_0$ induces a
graph of groups decomposition for $H$.

If $T_0$ is not locally finite, repeat.  As long as the resulting
Bass--Serre tree is not locally finite we can continue.  Since at each
step, we add one edge to the quotient graph of groups decomposition of
$G$, this process must terminate by Bestvina--Feighn \cite{ar:BF}.
\end{proof}

We also note that recently Kropholler has proved a more general
statement \cite{ar:K2}.


\section{Realization}\label{sc:rvts}

In this section we prove the main theorem.  For the remainder of this
paper, we let $G \in \LF$ be fixed and \GDG \ denote the
$\Out(G)$--invariant deformation space discussed in Section
\ref{sc:lf-groups}.

Suppose that $K$ is a subgroup of $\Out(G)$ and $K$ fixes some point
$T \in \GDG$.  Then $\tK$ the lift of $W$ to $\Aut(G)$ consists of
automorphisms $\phi$ such that there exists an isometry $h_\phi\co T \to
T$ where $h_\phi(gx) = \phi(g)h_\phi(x)$ for all $x \in T,g \in G$.
As the $G$--trees in \GDG \ are irreducible and minimal, $h_\phi$ is
unique \cite{ar:CM}.  Thus we get a homomorphism $\tK \to \Isom(T)$,
i.e. $T$ is a $\tK$--tree.  It is easy to check that $\tK$ extends the
action of $G/Z(G)$ on $T$.

As $T$ is a locally finite tree, if the edge groups are virtually
polycyclic then $\tK \in GVP(n')$ for some $n'$. In this case Theorem
\ref{th:st} implies that whenever \tK \ is a finite index subgroup of
some group \tH, then $\tH \in GVP(n')$.  Thus we have a $G$--tree
fixed by $H$, the image of \tH \ in $\Out(G)$.  We compute the edge
stabilizers for the action of \tK \ on $T$ via the following
sequences.  For an edge $f \subseteq T$ denote by $G_f$ (respectively
$\tK_f$) the edge stabilizer of $f$.

\begin{lemma}\label{lm:edge_stab}
The edge stabilizers of $T$ for the \tK --action, $\tK_f$, fit into
short exact sequences:
\begin{equation}\label{eq:ses_edge}
\xymatrix@R=0.0in{ 
&1 \ar[r] &G_f / Z(G) \ar[r] &\tK_f \ar[r] &K_f \ar[r] &1 \\}
\end{equation}

\noindent
where $W_f$ is the image of $\tK_f$ in $W.$ In particular, if $K_f$ is
virtually polycyclic then $\tK_f$ is a virtually polycyclic subgroup
of \tK.
\end{lemma}

\begin{proof}
This only place where exactness needs to be checked is that $G_f /
Z(G)$ is the kernel of the map $\tK_f \to K_f$.  This follows as $\tK$
extends the action of $G,$ hence $G / Z(G) \cap \tK_f = G_f / Z(G)$.
\end{proof}

We can now prove the main theorem.

\vskip 0.2 in
\noindent
\textbf{Main Theorem.} \textit{Suppose $G \in \LF$ and $K$ is a
polycyclic subgroup of $\Out(G)$ which fixes a point in \GDG.  If $H$
is a subgroup of $\Out(G)$ commensurable with $K$, then $H$ fixes a
point in \GDG.}

\begin{proof}
Suppose that $H \subseteq \Out(G)$ contains $K$ as a finite index
subgroup where $K$ is polycyclic and fixes a point in \GDG.  We then
have the following short exact sequences:
\begin{equation}\label{eq:ses_fs}
\xymatrix@R=0.0in{ 
&1 \ar[r] &G / Z(G) \ar[r] &\Aut(G) \ar[r] &\Out(G) \ar[r] &1 \\
&         &\Vert    &\bigvee        &\bigvee        & \\ 
&1 \ar[r] &G / Z(G) \ar[r] &\tH \ar[r] &H \ar[r] &1\\
&         &\Vert    &\bigvee       &\bigvee          & \\ 
&1 \ar[r] &G / Z(G) \ar[r] &\tK \ar[r] &K \ar[r] &1}
\end{equation}

\noindent
Then as $\tK \in GVP(n')$ by Lemma \ref{lm:edge_stab}, and as $\tH$
contains $\tK$ as a finite index subgroup, we have $\tH \in GVP(n')$
by Theorem \ref{th:st}.  Thus \tH \ acts on a locally finite tree $T'$
inducing an action of $G / Z(G)$, hence also $G$, on $T'$ with
virtually polycyclic stabilizers necessarily of Hirsch length $n$ by
Lemma \ref{lm:prop}.  Let $T$ be the minimal subtree of $T'$ for $G$.
Then $T \in \GDG$ and clearly $H$ fixes this $G$--tree.
\end{proof}

For $n=0$ or $1$, if $G \in \LF$ then for any point $T \in \GDG$ the
vertex and edge groups have finite outer automorphism groups.  Levitt
\cite{ar:L} has shown that for these $G$--trees the stabilizer is
virtually finitely generated abelian, hence we have the following
corollary:

\begin{corollary}\label{col:n=01}
If $G \in \LF$ for $n=0$ or 1 and $H$ is a subgroup of $\Out(G)$ which
contains a finite index subgroup that fixes a point in \GDG, then $H$
fixes a point in \GDG.
\end{corollary}


\providecommand{\bysame}{\leavevmode\hbox to3em{\hrulefill}\thinspace}


\end{document}